\documentclass[11pt]{article}
\usepackage{amssymb}
\usepackage{url}
\newcommand{\be}{\begin{equation}}
\newcommand{\ee}{\end{equation}}
\newtheorem{lemma}{Lemma}
\begin{document}
\title{\textbf{The Discrete Logarithm Problem over Prime Fields can be transformed to a Linear Multivariable Chinese Remainder Theorem}}
\author{H. Gopalakrishna Gadiyar and R. Padma \\
Department of Mathematics\\
School of Advanced Sciences\\ V.I.T. University, Vellore 632014 INDIA\\E-mail: gadiyar@vit.ac.in, rpadma@vit.ac.in}
\date{~~}
\maketitle

\begin{abstract}
We show that the classical discrete logarithm problem over prime fields can be reduced to that of solving a system of linear modular equations.


\end{abstract}

\section{Introduction} The first published public key cryptographic algorithm is the famous Diffie - Hellman key exchange protocol which is based on the intractability of solving the discrete logarithm problem over prime fields for large primes. The discrete logarithm over prime fields is defined as follows: Let $p>2$ be a prime and $a_0$ be a primitive root of $p$. We know that every $b_0 \in \{1, 2, \cdots ,p-1\}$ can be expressed as a power of $a_0$ mod $p$. That is,  
\be 
a_0^n \equiv b_0 {\rm{~mod~}}p \label{eq:basiceqn}
\ee 
for a unique $n$ modulo $p-1$. Then $n$ is called the discrete logarithm or index of $b_0$ with respect to the base $a_0$ modulo $p$. Finding $n$ modulo $p-1$ given $a_0$ and $b_0$ modulo $p$ is called the discrete logarithm problem over prime fields. If $p$ is a randomly chosen large prime, it is believed that this problem is computationally infeasible and hence is used as the basis of the Diffie - Hellman key exchange protocol. 

\section{Earlier paper} In an investigation carried out earlier \cite{gadiyar:maini:padma}, the authors had derived a linear modular equation in two unknowns. From (\ref{eq:basiceqn}), we get
\be
a_0^{np} \equiv b_0^p {\rm{~mod~}}p^2 \, . 
\ee
This can be written as
\be
(a_0+a_1p)^n \equiv b_0+b_1p {\rm{~mod~}}p^2 \, . \label{eq:teich}
\ee
Note that here $a_0+a_1p$ and $b_0+b_1p$ are the truncation of the Teichmuller expansions got by Hensel lifting the polynomial equation $x^p-x=0$. Using the binomial theorem, this equation gets linearized. So (\ref{eq:teich}) becomes
\be
a_0^n + na_0^{n-1}a_1p \equiv  b_0+b_1p {\rm{~mod~}}p^2 \, .
\ee
Now introducing a new term $\beta_n$ which is the ``carry" 
\be
a_0^n - b_0 \equiv \beta_n p {\rm{~mod~}}p^2  \, ,
\ee
we get a linear equation in two variables
\be
\beta_n +n \frac{b_0}{a_0}a_1 \equiv b_1 {\rm{~mod~}}p \, .
\ee 
Thus if we know $a_0^n$ mod $p^2$ where $1 \le n \le p-1$, then we can find $n$. In this connection, see \cite{Catalano:Nguyen:Stern}. While $\beta _n$ is a multiplicative carry, the additive carries and their connections to many areas of mathematics are studied in \cite{Diaconis}, \cite{Diaconis:Fulman}, \cite{Hanlon} and \cite{Holte}.

We now modify the argument because the equation (\ref{eq:basiceqn}) for the discrete logarithm is given modulo $p$ and the index $n$ has to be found modulo $p-1$. Hence we study the discrete logarithm problem for composite modulus.

\section{The New Discrete Logarithm Problem} To make the ideas available to a larger audience we take the safe prime case. $p$ is a safe prime if $p=2q+1$ where $q$ is a prime. Then we have the following lemma.

\begin{lemma} Let $gcd(a_0,q)=1$ and $gcd(b_0,q)=1$. Let $a_0$ be a primitive root of $p$ and let $a_0$ and $b_0$ satisfy (\ref{eq:basiceqn}). Then
\be
a_0^{n\phi(q)} \equiv b_0^{\phi(q)} {\rm{~mod~}}pq \, .
\ee

\end{lemma}

\noindent{\bf Proof:} Raising both sides of (\ref{eq:basiceqn}) to the power of $\phi(q)$ one gets
\be
a_0^{n\phi(q)} \equiv b_0^{\phi(q)} {\rm{~mod~}}p \, ,
\ee
and by Fermat's little theorem
\be
a_0^{n\phi(q)} \equiv b_0^{\phi(q)} ~(\equiv ~1){\rm{~mod~}}q \, ,
\ee
trivially. Since $(p,q)=1$, the lemma follows.

Note that the subgroup generated by $a_0$ mod $pq$ is of order $q$ and hence we can find $n$ modulo $q$. Thus $n$ modulo $p-1$ is $n$ or $n+q$ modulo $p-1$. 

In \cite{lerch}, Lerch defined the Fermat quotient for a composite modulus. Let $x$ be such that $gcd(x,n)=1.$ Then $q(x)$ defined by
\be
x^{\phi(n)} \equiv 1+q(x)n {\rm{~mod~}}n^2 \, .
\ee
is called the Fermat quotient of $x$ modulo $n$. We do not use the Euler's $\phi $-function but we use Carmichael's $\lambda $ function. $\lambda (n)$ is defined as follows \cite{Cameron:Preece}. $\lambda(2)=1$, $\lambda(4)=2$ and
\be
\lambda(n)=~\left \{ \begin{array}{ll} \phi(p^r),~&{\rm ~if ~} n=p^r\\
                                          2^{r-2},&{\rm ~if ~} n=2^r, ~r\ge 3\\
																				lcm(\lambda(p_1^{r_1}), \lambda(p_2^{r_2}),\cdots , \lambda(p_k^{r_k})),& {\rm ~if ~} n=p_1^{r_1}p_2^{r_2}\cdots p_k^{r_k} \end{array} \right .
\end{equation}
When $n=p^2q^2$, where $p=2q+1$, $q$ is a prime, $\phi(p^2q^2)=2pq^2\phi(q)$ and $\lambda(p^2q^2)=pq\phi(q)$. In other words, the order of the group of units modulo $p^2q^2$ is $\phi(p^2q^2)$, whereas the order of the largest cyclic group modulo $p^2q^2$ is $\lambda(p^2q^2)$. Hence we define $q(x)$ by the congruence
\be
x^{pq\phi(q)} \equiv 1+q(x)p^2q^2 {\rm{~mod~}}p^3q^3 \, .
\ee
Now we have the analogue of Teichmuller expansion (\ref{eq:teich}) modulo $p^2q^2$.  

\begin{lemma}  Let $(x)_l$ denote the residue of $(x)$ modulo $l$. Then,
\be
(a_0^{\phi(q)}+a_1pq)^n \equiv (b_0^{\phi(q)})_{pq}+b_1pq) {\rm{~mod~}}p^2q^2 \, , \label{eq:Eqn}
\ee
where $a_1=-a_0^{\phi(q)}q(a_0)$ mod $pq$ and $b_1=-b_0^{\phi(q)}q(b_0)$ mod $pq$.
\end{lemma}

\noindent{\bf Proof} We want $a_1$ and $b_1$ to satisfy (\ref{eq:Eqn}). Using the carry notation 
\be
a_0^{n\phi(q)} \equiv (b_0^{\phi(q)})_{pq}+\beta_n pq {\rm{~mod~}}p^2q^2 \, , \label{eq:neweQ}
\ee
we get the equation
\be
\beta_n +n \frac{b_0^{\phi(q)}}{a_0^{\phi(q)}} a_1 \equiv b_1 {\rm{~mod~}}pq \, . \label{eq:neweq}
\ee
Taking $pq^{th}$ power on both sides of (\ref{eq:neweQ}),
\be
a_0^{npq\phi(q)} \equiv ((b_0^{\phi(q)})_{pq}+\beta_n pq)^{pq} {\rm{~mod~}}p^3q^3
\ee
and using the definition of $q(x)$, we get
\be
nq(a_0) \equiv q(b_0)+\frac{\beta_n}{b_0^{\phi(q)}} {\rm{~mod~}} pq \,. \label{eq:neweq1}
\ee
Comparing (\ref{eq:neweq}) and (\ref{eq:neweq1}) will give the desired values of $a_1$ and $b_1$.

\noindent {\bf Remark:} Note that $a_1$ and $b_1$ in (\ref{eq:Eqn}) can be calculated in polynomial time.

Taking equation (\ref{eq:neweq}) mod $p$ and $q$, we get
\begin{eqnarray}
\beta_n +n \frac{b_0^{\phi(q)}}{a_0^{\phi(q)}}a_1 &\equiv &b_1 {\rm{~mod~}}p \\
\beta_n+na_1&\equiv &b_1 {\rm{~mod~}}q \, .
\end{eqnarray}
Hence the discrete logarithm problem can be transformed to the multivariable Chinese remainder theorem. 

\section{Numerical example} 
To make the numerical work easy and understandable, we take the prime $p=11$ and $q=5$. $2$ is a primitive root of $11$. We take $a_0=2$ and $b_0=4$. We get $q(a_0)=18$ and $q(b_0)=36$, $a_1=42$ and $b_1=28$. We get the linear congruence
\be
\beta_n+12n \equiv 28 {\rm{~mod~}}55 \, ,
\ee
which becomes two simultaneous linear congruences in two unknowns with relatively prime moduli.
\begin{eqnarray}
\beta_n+n &\equiv & 6 {\rm{~mod~}} 11 \, , \\
\beta_n+2n &\equiv &3 {\rm{~mod~}}5 \, .
\end{eqnarray}

\section{Historical comments and Conclusion} The inspiration for this paper is the successful attack of elliptic curve discrete logarithm problem on anomalous elliptic curves \cite{sat:ara}, \cite{semaev} and \cite{smart}. For a discussion of the discrete logarithm problem modulo a composite integer, see \cite{bach}. A formula for solving the discrete logarithm problem in certain cases was obtained by Riesel \cite{riesel} using Fermat quotient and its generalisations. 

The authors do not possess a library of classical books on number theory. We therefore refer to the work of Professor Oliver Knill \cite{Knill} of Harvard. He says that the multivariable Chinese remainder theorem has not been investigated thoroughly. 

It is a pleasant surprise to the authors that the fundamental problem of public key cryptography which started with the Diffie - Hellman key exchange protocol has been connected to a classical result with a hoary past.

This paper is dedicated to S. Ramanathan on his birth centenary. The second author is his daughter and the first author is his son-in-law.

\end{document}